\documentclass[leqno,12pt]{article} 
\usepackage{amsmath, amssymb}
\usepackage{amsthm} 
\usepackage[all]{xy}
\usepackage{amssymb}
\usepackage{graphicx}
\usepackage{amsbsy}
\usepackage{amsfonts}
\usepackage{mathrsfs}
\usepackage{amstext}
\usepackage{amscd}
\usepackage[dvips]{epsfig}
\usepackage{psfrag}
\usepackage{enumerate}
\usepackage{flafter}
\usepackage{color}
\usepackage{url}
\usepackage{here}
\allowdisplaybreaks

\newcommand{\C}{{\Bbb C}}
\newcommand{\Q}{{\Bbb Q}}

\newcommand{\RR}{\mathcal{R}}

\newcommand{\lra}{\longrightarrow}

\newcommand{\Z}{{\Bbb Z}}
\newcommand{\om}{\omega}
\newcommand{\x}{\mathrm{x}}
\theoremstyle{plain} 
\newtheorem{theorem}{\indent\sc Theorem}[section] 
\newtheorem{lemma}[theorem]{\indent\sc Lemma}

\theoremstyle{definition} 
\newtheorem{definition}[theorem]{\indent\sc Definition}

\newtheorem{example}[theorem]{\indent\sc Example}

\newtheorem{problem}[theorem]{\indent\sc Problem}
\makeatletter
\def\address#1#2{\begingroup
\noindent\parbox[t]{7.8cm}{%
\small{\scshape\ignorespaces#1}\par\vskip1ex
\noindent\small{\itshape E-mail address}%
\/: #2\par\vskip4ex}\hfill%
\endgroup}%
\makeatother
\title{\uppercase{Relationship between quandle shadow cocycle invariants and Vassiliev invariants}} 
\author{
\textsc{Sukuse Abe$^{*}$} 
}
\date{} 
\begin{document}

\maketitle


\newsavebox{\boxa}
\sbox{\boxa}{\includegraphics[height=2cm]{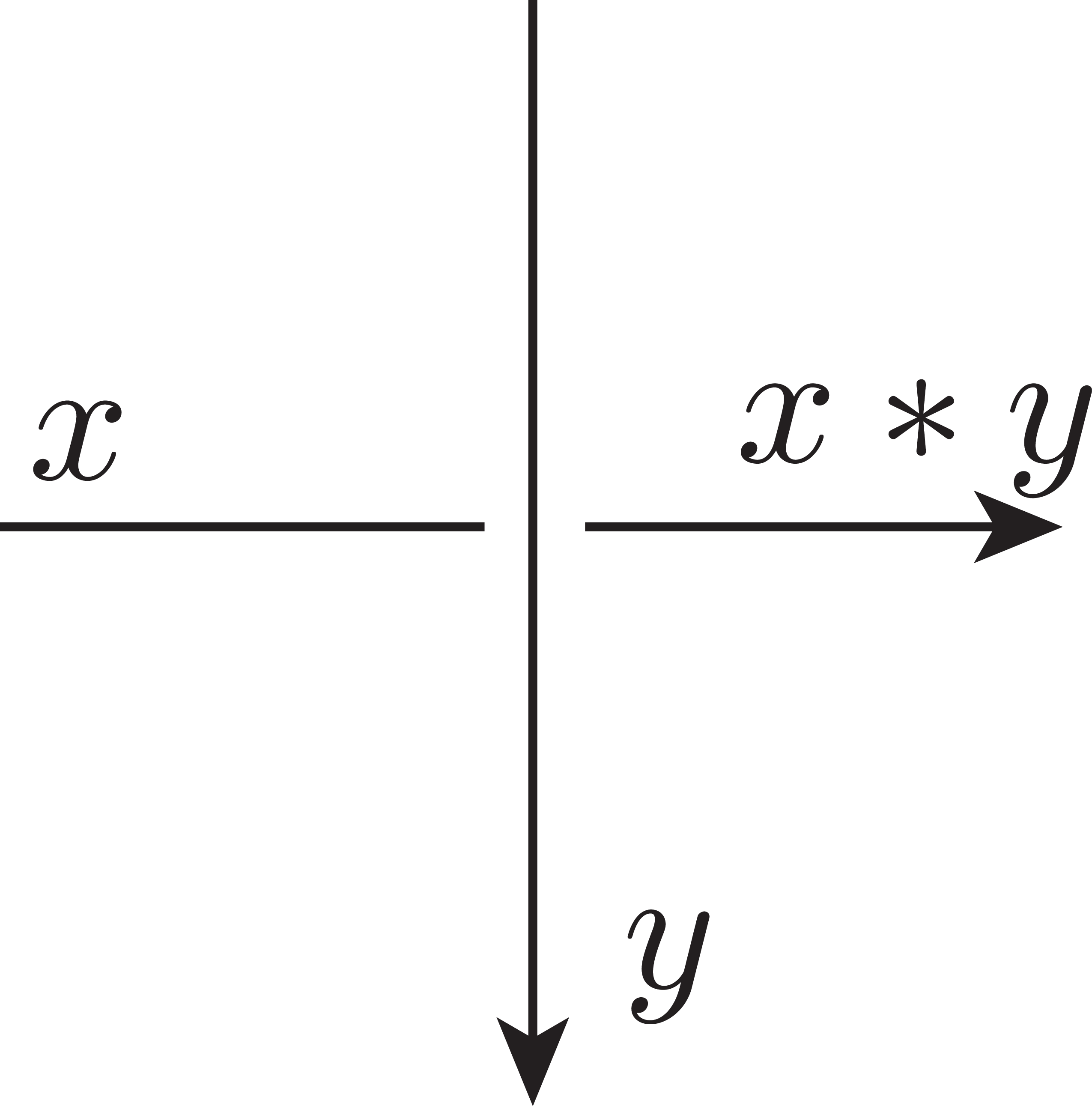}}
\newlength{\bw}
\settowidth{\bw}{\usebox{\boxa}}
\newsavebox{\boxb}
\sbox{\boxb}{\includegraphics[height=2cm]{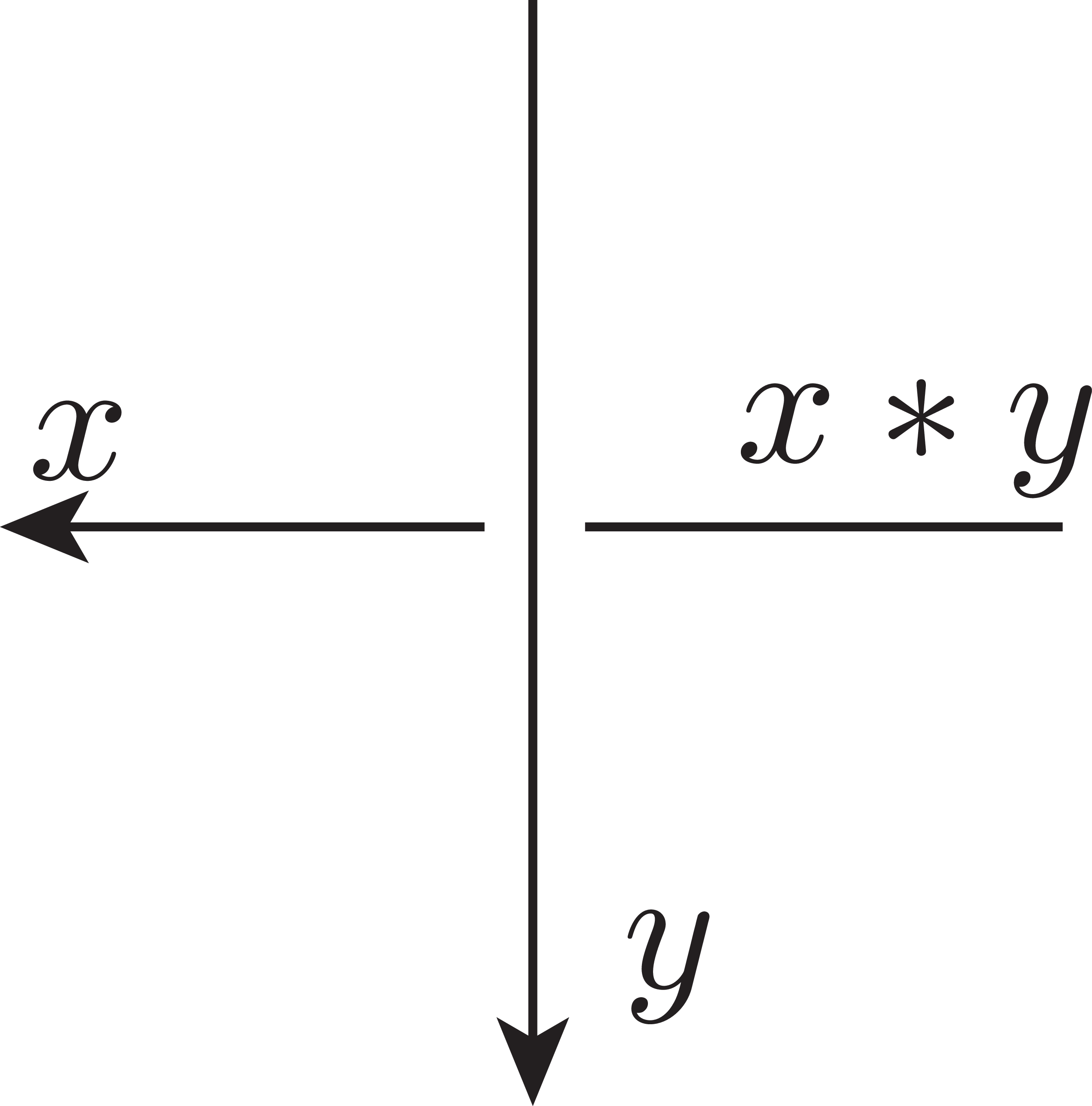}}
\settowidth{\bw}{\usebox{\boxb}}
\newsavebox{\boxc}
\sbox{\boxc}{\includegraphics[height=2cm]{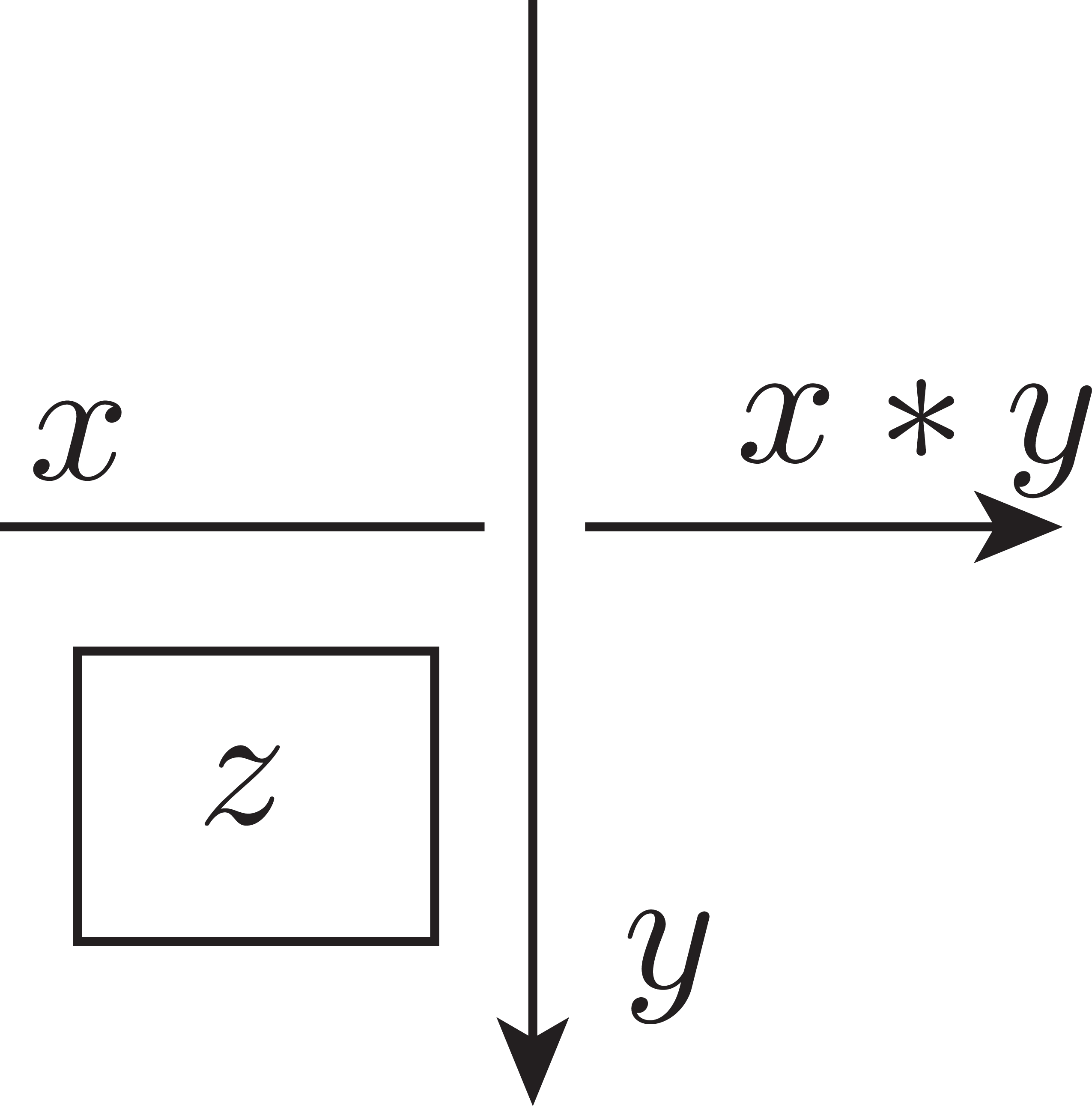}}
\settowidth{\bw}{\usebox{\boxc}}
\newsavebox{\boxd}
\sbox{\boxd}{\includegraphics[height=2cm]{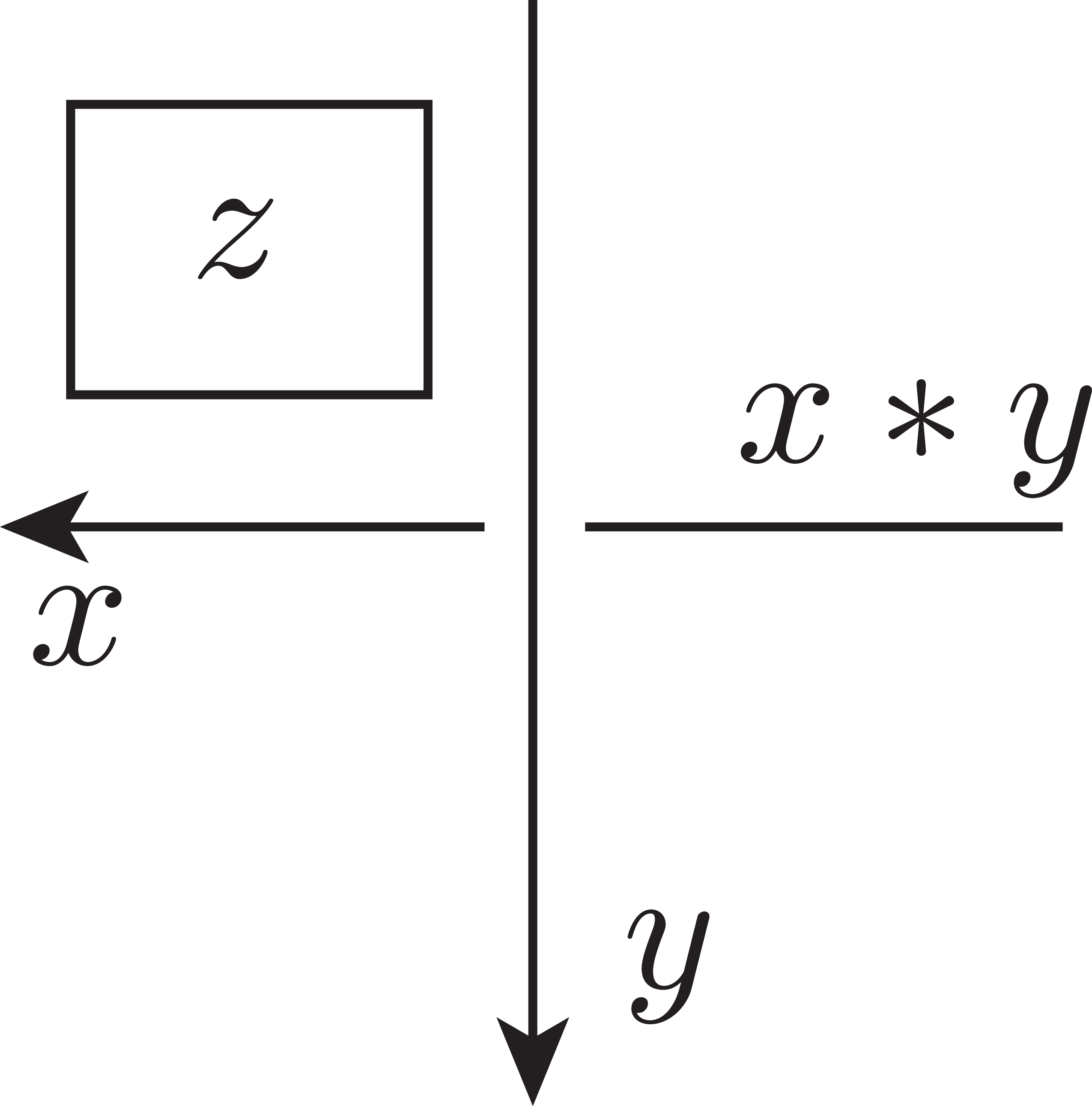}}
\settowidth{\bw}{\usebox{\boxd}}

\begin{abstract}
As one of the problems in his list \cite{16:2002},
T. Ohtsuki proposed to study relations between quandle cocycle invariants and quantum invariants.
The aim of this paper is to answer one of those questions.
We prove that the coefficient of the finite perturbative expansion of the quandle shadow cocycle invariant
defined by $(\Z /p\Z)$-laurent polynomial quandle is Vassiliev invariant for any links.
\end{abstract}

\section{Introduction} 
\label{sec:introduction}
The present work is motivated by \cite[\S 5.4 Quandle cocycle invariants]{16:2002},
in particular the following problem:
\begin{problem}[{\cite[Problem 5.7]{16:2002}}]
Find relations between quandle cocycle invariants and knots invariants known so far, such as quantum invariants.
\end{problem}
{M. Gra\~{n}a proved that the quandle cocycle invariants can be presented by knot invariants derived from certain ribbon categories \cite{24:2002}.
A. Soloviev, develop a theory of non-unitary set-theoretical solutions to the
Quantum Yang--Baxter equation \cite{18:2000}. 
S. Abe obtained Vassiliev invariants in the case where quandle (shadow) cocycles are calculated by using trivial quandles \cite{1:2016}.
We knew only these three results, and Vassiliev invariants and quandle (shadow) cocycle invariants were thought to have little relation to each other.

We obtain the following main Theorem:
\begin{theorem}\label{100}
Let $p$ be an odd prime,
let $\xi={\rm exp}(2 \pi \sqrt{-1}/p)$,
let $\Z/p\Z \cong \langle t \mid t^p=1 \rangle$,
let $\phi'$ be a $3$-cocycle of a Alexander quandle $(\Z/p \Z, \om)$,
and let $\Phi_{\phi'}(L, \om)\in \Z[\Z /p\Z]\cong \Z[t]/(t^p-1)$
be a quandle shadow cocycle invariant of a link $L$.
A quandle shadow cocycle invariant $\Phi_{\phi'}(L,\om)$ has the exponent a Lauren polynomial
$\Z[\om,\om^{-1}]$ of $t$.
We have that
\begin{align*}
&\Phi_{\phi'}(L,\om)|_{t=\xi}\\
&\equiv \sum_{d=0}^{p-2}a_{p,d}(L,\om)(\xi-1)^d \pmod{(\xi-1)^{p-1}}\\
&\in \Z[\xi,\frac{1}{2},\frac{1}{3},\ldots ,\frac{1}{p-2},\phi']\\
&\cong \Z[t,\frac{1}{2},\frac{1}{3},\ldots ,\frac{1}{p-2}, \phi']/(1+t+t^2+\cdots +t^{p-1}).
\end{align*}
When we substitute 
$1/(1-\hbar)$ for $\om$ of $a_{p,d}(L,\om)\in \Z/p\Z[\om,\om^{-1}] \ (0\leq d\leq p-2)$,
the coefficient of $\hbar^k$ in
$
\displaystyle a_{p,d}(L,\om)|_{Y=1/(1-\hbar)}\in \Z/ p \Z[[\hbar ]]
$
is a $(\Z/p\Z)$-valued Vassiliev invariant of degree $k$ for any $0\leq d\leq p-2$.
\end{theorem}

For any braids, we can obtain Vassiliev invariants from quandle shadow cocycle invariants. 
This is expected to have applications to surface links and low-dimensional manifolds in the future.

In Section \ref{200}, we review the quandle (shadow) cocycle invariants.
A {\it quandle} is a set with a binary operation satisfying certain axioms,
which is analogous to a group with conjugation.
The cohomology groups of quandles were introduced by 
Carter, Jelsovsky, Kamada, Langford and Saito \cite{3:2003}
as an analogy of group cohomology.
The $n$-th cohomology group 
of quandle $X$ with coefficient group $A$ shall be denoted by $H^n_Q(X;A)$.
Furthermore, the quandle shadow cocycle invariants of classical links 
were defined using $3$-cocycles of the cohomology groups of quandles \cite{3:2003, 4:2001}.

In Section \ref{201}, we review the Vassiliev invariants.
The notion of {\it Vassiliev invariant} was introduce by V.A. Vassiliev \cite{20:1990}.
We know that after a suitable change of variables each coefficient in the Taylor expansion of
the Jones polynomial is a Vassiliev invariant \cite{25:1993}.

In Section \ref{202}, we drive finite perturbative expansion from quandle shadow cocycle invariants.
It is the first time that we have ever derived 
finite perturbative expansion from quandle shadow cocycle invariants.

In Section \ref{203}, we prove that the Vassiliev invariants can be deduced quandle shadow cocycle invariants by using by the method similar to how
Vassiliev invariants are derived from Quantum polynomial invariants.
A new $(\Z/p \Z)$-valued Vassiliev invariant is derived from one value of 
quandle shadow cocycle invariant by using an odd prime $p$.

\section{Quandle and Quandle shadow cocycle invariants}\label{200}

Firstly,
we review the definitions of quandles and their cohomology groups of 
S. Carter, D. Jelsovsky, S. Kamada, L. Langford, M. Saito $\cite{3:2003}$.

A {\it quandle} is a set with a binary operation satisfying certain axioms,
which is analogous to a group with the conjugation.
It is a set with a self-distributive binary operation whose definition was motivated
from knot theory.

A quandle is a set $X$ with a binary operation 
$\ast:X\times X\to X$ such that the following three conditions are satisfied:
\begin{enumerate}[{(i)}]
\item For any $a\in X$, $a\ast a=a.$
\item For any $a,b \in X$, there exists a unique $c \in X$ 
such that $a=c\ast b$.
\item For any $a,b,c \in X$, $(a\ast b)\ast c=(a\ast c)\ast(b\ast c).$
\end{enumerate}

Cohomology groups of quandles were introduced by 
S. Carter, D. Jelsovsky, S. Kamada, L. Langford, M. Saito $\cite{3:2003}$,
as an analogy of group cohomology;
we denote by $H^n_Q(X;A)$ the $n$-th cohomology group 
of a quandle $X$ with a coefficient group $A$.

We review the quandle cochain complex \cite{3:2003}. 
Let $(X,\ast)$ be a finite quandle and $A$ be an abelian group. Consider an abelian group,
\[C_Q^n(X;A):=\{ \ f:X^n\to A \mid \ f(x_1,\dots,x_n)=0 \ {\rm when} \ x_i=x_{i+1} \ 
{\rm for} \ {\rm some} \, \, i\} \ .\]
For $n\geq 1$, we define the coboundary map of the set above, $\delta_n:C_Q^n(X;A)\lra C_Q^{n+1}(X;A)$, as follows:\\
\begin{align*}
\delta_n(f)(x_1,\dots,x_{n+1}):=&\sum_{i=2}^{n+1}(-1)^i\bigl(f(x_1,\dots,x_{i-1},x_{i+1},\dots,x_{n+1})\\
&-f(x_1\ast x_i,\dots,x_{i-1}\ast x_i,x_{i+1},\dots,x_{n+1})\bigr).
\end{align*}
We easily see that $\delta_{n+1}\circ\delta_n=0$. 
The cohomology of this complex $(C_Q^n(X;A),\delta_n)$ is denoted by $H_Q^n(X;A)$ and called the 
{\it quandle cohomology} of $X$ with coefficient group $A$.

%
%
%

Further, S. Carter, D. Jelsovsky, S. Kamada, L. Langford, M. Saito defined invariants of links
using the quandle cocycle.
Especially, by using $3$-cocycles of the cohomology groups of quandles,
{\it quandle shadow cocycle invariants} of links can be defined
\cite{3:2003, 4:2001}.
To define the quandle shadow cocycle invariants, we first need to define the {\it $X$-shadow coloring}.

Let $D$ be a diagram of the oriented link $L$,
$\x$ be a crossing of the diagram $D$, and $\mathcal{A}(D)$ be the set of arcs of $D$.
An {\it $X$-coloring} \cite[Definition.4.1]{4:2001} of $D$ is a map $C:\mathcal{A}(D) \to X$ satisfying the   
condition $C(\gamma_1)\ast C(\gamma_2)=C(\gamma_3)$ in  Figure \ref{10002} at each crossing $\x$
of $D$.
If $X$ is a finite quandle, then the number of $X$-colorings of $D$ is an invariant of link $L$ \cite{7:1982}.
\begin{figure}[htb]
	\centering
    \includegraphics[height=2.5cm]{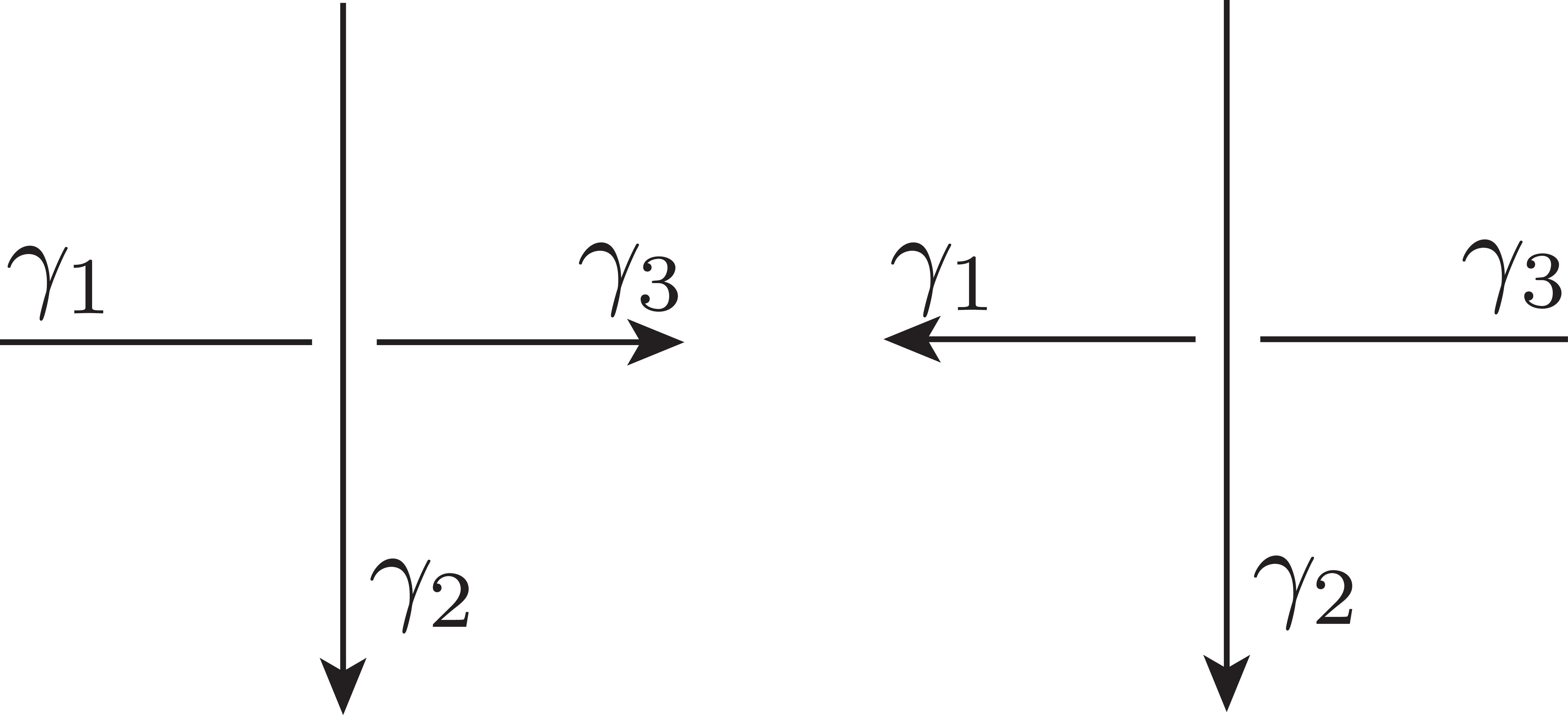}
    \caption{The coloring at each crossing $\x$.}
    \label{10002}
\end{figure}
Let $R(D)$ be the set of regions of the underlying immersed curves of $\mathcal{A}(D)$. 
The map $C:\mathcal{A}(D) \sqcup R(D) \to X$ is an $X$-shadow coloring of $D$ [{Definition.4.3 \cite{4:2001}}], as shown in 
Figure \ref{10002} for $C|_{\mathcal{A}(D)}:\mathcal{A}(D)\to X$ and in Figure \ref{10005} for $C|_{R(D)}:R(D)\to X$.
\begin{figure}[H]
	\begin{center}
    \includegraphics[height=2.5cm]{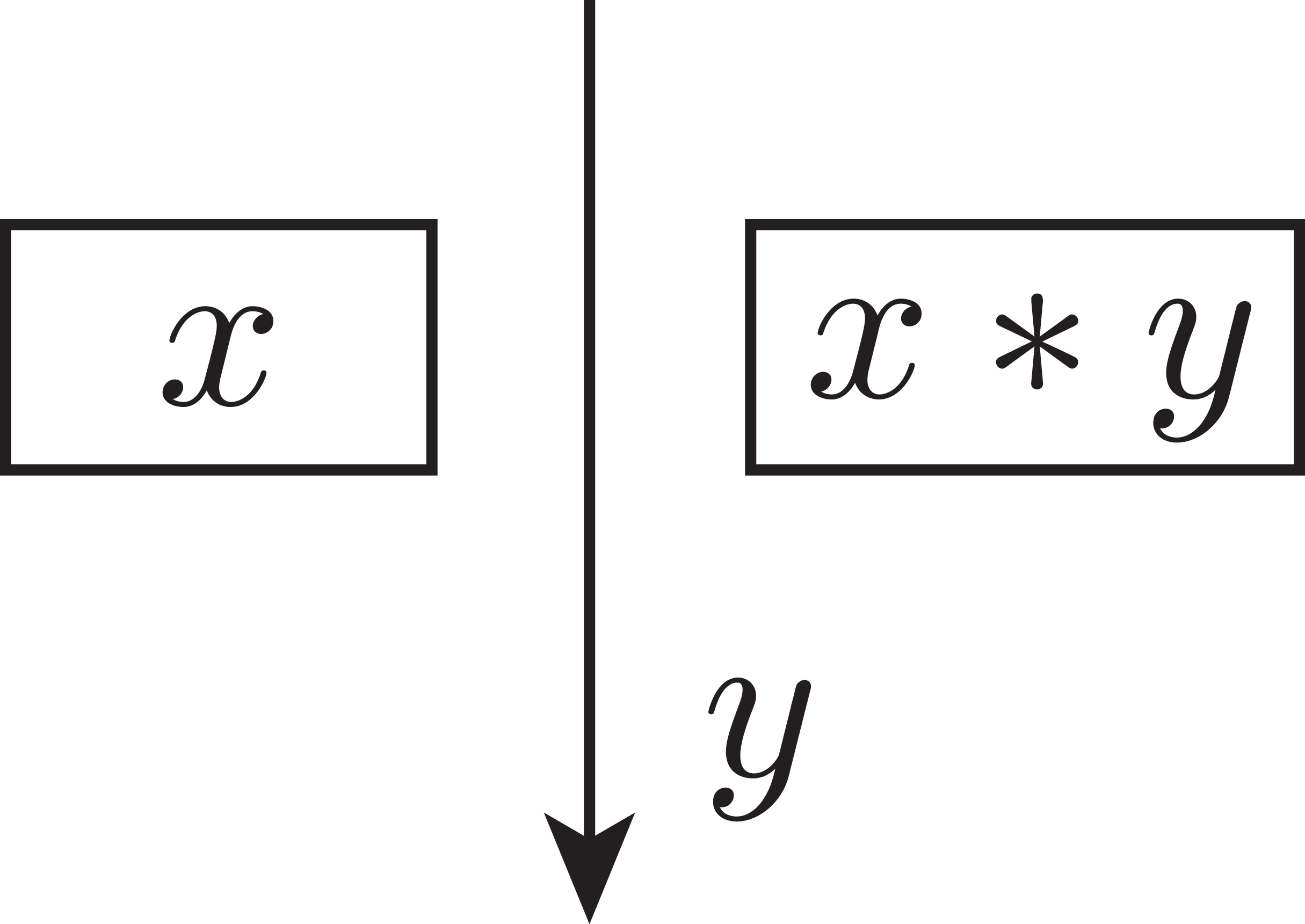}
    \caption{Boxes indicate the colorings of regions}
    \label{10005}
    \end{center}
\end{figure}

Let a finite quandle $X$ and a map $\phi:X^3\to A$ be given.
A ({\it Boltzmann}) {\it weigh} $W_{\phi}(\x; C)$, at crossing $\x$ is defined as follows:
$$
W_{\phi}\Biggl( \parbox{\bw}{\usebox{\boxc}} \Biggr)=\phi(z,x,y)\in A,
$$
$$
W_{\phi}\Biggl( \parbox{\bw}{\usebox{\boxd}} \Biggr)=\phi(z,x,y)^{-1}\in A.
$$
The following state sum is: 
$$
\sum_{C}\prod_{\x}W_{\phi}(\x; C)\in \Z[A].
$$
If $\phi:X^3\to A$ is a quandle $3$-cocycle of $X$, then the above state sum 
is an invariant of link $L$. We denote this invariant by
$$
\Phi_{\phi}(L).
$$
$\Phi_{\phi}(L)$ is called the {\it quandle shadow cocycle invariant} \cite{22:2004}.

S. Asami and S. Satoh calculated the quandle shadow cocycle invariant of torus links \cite{2:2005}.
M. Iwakiri calculated the quandle shadow cocycle invariant of two bridge links \cite{I}.

\section{Vassiliev invariants}\label{201}

Secondly, we explain the definition of Vassiliev invariants.
Let $\mathcal{K}$ be a vector space over $\C$ freely spanned by the isotopy classes of oriented knots in $S^3$. 
A {\it singular knot} is an immersion of $S^1$ into $S^3$, whose singularities are transversal double points. 
We regard a singular knot as a linear sum in $\mathcal{K}$ obtained by the relation shown in the following Figure \ref{10006}.
\begin{figure}[htb]
 \centering
 \includegraphics[width=5.5cm]{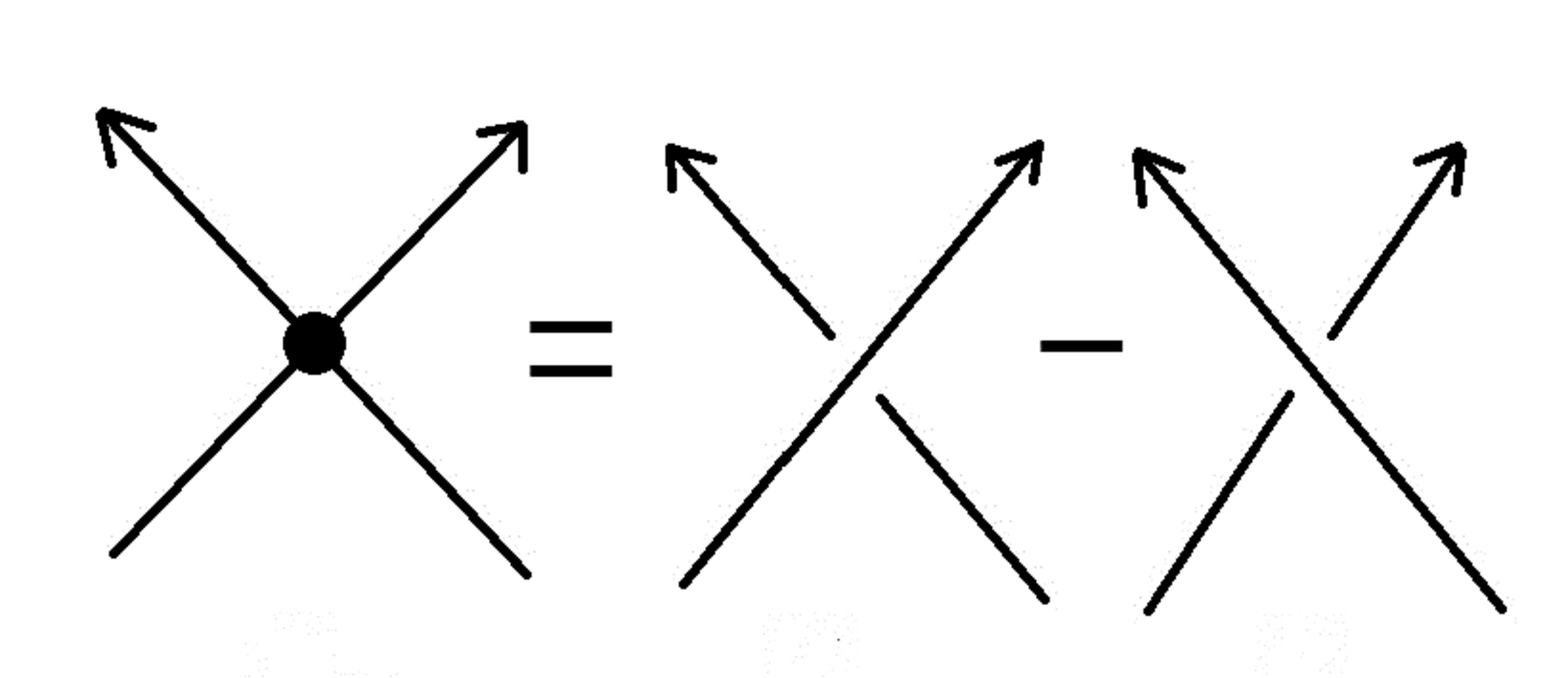}
 \caption{\rm Singular point.} \label{10006}
\end{figure}

Let $B$ be an abelian group,
$\mathcal{K}_d$ denote the vector subspace of $\mathcal{K}$ spanned by singular knots with $d$ double points. 
A homomorphism map $v:\mathcal{K}\to B$ is called a {\it Vassiliev invariant} of {degree $d$} if 
$v|_{\mathcal{K}_{d+1}}=0$ \cite{5:1994,6:1994,20:1990}.

The construction of the general quantum invariants is explained below.
\begin{theorem}[{\cite{2:1947}}]
The braid group on $n$ strands can be presented as
$$
B_n:=\langle \sigma_1,\ldots ,\sigma_{n-1} \mid \sigma_i\sigma_j=\sigma_j\sigma_i \  \text{for} \ 
|i-j|\geq 2, \ \ \sigma_i\sigma_j\sigma_i=\sigma_j\sigma_i\sigma_j \  \text{for}  \ |i-j|=1\rangle,
$$
where the braid $\sigma_i$ performs a positive crossing of the strands $i$ and $i+1$.
Similarly, the braid $\sigma_i^{-1}$ performs a negative crossing of the strands $i$ and $i+1$.
\end{theorem}
\begin{definition}[{\cite{9:1989}}]
Let $V$ be a vector space over $\C$.
We obtain a representation, $\psi_n:B_n \to {\rm End}(V^{\otimes n})$, defined as
\begin{equation}\label{10008}
\psi_n(\sigma_i)=(id_V)^{\otimes (i-1)}\otimes R \otimes (id_V)^{\otimes(n-i-1)}.
\end{equation}
Such a map $\psi_n$ given in (\ref{10008}) always satisfies $\psi_n(\sigma_i \sigma_j)=\psi_n(\sigma_j \sigma_i) \, \, (|i-j|\geq 2)$.
To obtain $\psi_n(\sigma_i\sigma_{i+1}\sigma_i)=\psi_n(\sigma_{i+1}\sigma_i \sigma_{i+1})$, 
the matrix $R$ is required to satisfy the following equation:
$$
(R\otimes id_V)(id_V \otimes R)(R \otimes id_V)=(id_V \otimes R)(R \otimes id_V)(id_V \otimes R ).
$$
We call this equation the {\it Yang-Baxter equation}, and its solution is called an {\it $R$-matrix}.
For an $R$-matrix we obtain a representation $\psi_n$ of the braid group $B_n$ given by (\ref{10008}).
\end{definition}
\begin{theorem}[{chap.I \cite{19:1994}} and {chap.X \cite{10:1995}}]\label{10009}
Let $L$ be an oriented link and $b\in B_n$ be a braid such that the closure is isotopic to $L$.
If an R-matrix $\RR$ and a linear map $h\in {\rm End}(V)$
satisfy the following relations:
\begin{itemize}
\item ${\rm trace}_2\bigl((id_V \otimes h)\cdot \RR^{\pm}\bigr)=id_V$,
\item $\RR\cdot(h\otimes h)=(h\otimes h)\cdot \RR$,
\end{itemize}
then
${\rm trace}\bigl(h^{\otimes n}\cdot \psi_n(b)\bigr)$
is invariants of $L$.  
\end{theorem}
Jones polynomial, which is a quantum invariant, can be derived from the above invariant.

If
$$
\RR=\begin{pmatrix}
t^{1/2} & 0 & 0 & 0 \\
0 & 0 & t & 0 \\
0 & t & t^{1/2}-t^{3/2} & 0 \\
0 &0 & 0 & t^{1/2}\\
\end{pmatrix}
$$
and
$$
h=\begin{pmatrix}
t^{-1/2} & 0 \\
 0 & t^{1/2}\\
\end{pmatrix},
$$
then ${\rm trace}\bigl(h^{\otimes n}\cdot \psi_n(b)\bigr)$ is a Jones polynomial.
Moreover, Vassiliev invariants are derived from Jones polynomial by the following theorem.

\begin{theorem}[{\cite{25:1993}}]\label{15}
For any link $L$, let $J_L(t)$ be the Jones polynomial of $L$.
Then, the coefficient of $\hbar^k$ in $J_L(t)|_{t=e^{\hbar}}$ is a Vassiliev
invariant of degree $k$.
\end{theorem}
\begin{proof}
Let $D$ be a diagram of $L$. 
We associate matrix $R$ and its inverse, $R^{-1}$, to the positive and negative crossings of $D$ respectively.
We substitute $ t=e^{\hbar}$ in $J_L(t)$.
Since
$$
\RR=\begin{pmatrix}
t^{1/2} & 0 & 0 & 0 \\
0 & 0 & t & 0 \\
0 & t & t^{1/2}-t^{3/2} & 0 \\
0 &0 & 0 & t^{1/2}\\
\end{pmatrix}
$$
and
$$
\RR^{-1}=\begin{pmatrix}
t^{-1/2} & 0 & 0 & 0 \\
0 & t^{-1/2}-t^{-3/2} & t^{-1} & 0 \\
0 & t^{-1} & 0 & 0 \\
0 &0 & 0 & t^{-1/2}\\
\end{pmatrix},
$$
$\RR=\RR^{-1}$ when $\hbar=0$.
Therefore, $\RR-\RR^{-1}$ is a matrix whose entries are divisible by $\hbar$. 
This difference is associated to the double point of a singular link that occurs in the definition of Vassiliev invariants. 
Theorem \ref{10009} shows that $J_L(t)$ is composed by the product of R-matrix $\RR$.
Therefore, if $L$ is a singular link with exact $k+1$ singular points, 
then $J_{L}(t) |_{t=e^{\hbar}}$ is divisible by $\hbar^{k+1}$.
Hence, the coefficient of $\hbar^k$ is equal to $0$ for such singular links.
\end{proof}


\section{Finite perturbative expansion of quandle shadow cocycle invariants}\label{202}

Thirdly, quandle shadow cocycle invariants have finite perturbative expansion.
T. Le, J. Murakami and T. Ohtsuki showed that the quantum invariants of $3$-manifolds
have perturbative expansions and Fermat-limits \cite{27:1998}.
We review Isomorphism theorem necessary to show that quandle shadow cocycle invariants $\Phi_{\phi}(L)$
have finite perturbative expansions.

We prove that all coefficient are unique when $\Phi_{\phi}(L)|_{t=\xi}$ has a finite perturbative expansion.

%

Let $p$ be an odd prime and let $\xi={\rm exp}(2 \pi \sqrt{-1}/p)$.
We define the map
$$
\varphi: \Z[t] \to \Z[\xi]
$$
by $\varphi(f(t)):= f(\xi)$ i.e. just replacing $t$ by $\xi$.
Since $\varphi$ is a surjective homomorphism and
${\rm Ker} \ \varphi=(1+t+t^2+\cdots +t^{p-1})$,
we obtain the ring $\Z[t]/(T(t))$ is isomorphic to $\Z[\xi]$ where
$$
T(t)=1+t+t^2+\cdots +t^{p-1}.
$$
Since $\Phi_{\phi}(L)$ is an invariant of $L$,
$\Phi_{\phi}(L)|_{t=\xi}$ is also an invariant of $L$.

\begin{lemma}[{\cite{23:2002}, Lemma 9.7.}]\label{13}
Let $p$ be an odd prime and let $\xi={\rm exp}(2 \pi \sqrt{-1}/p)$. 
We consider $\Phi_{\phi}(L)|_{t=\xi} \in \Z[\xi]$. We put 
$$
\Phi_{\phi}(L)|_{t=\xi}=a_{p,0}(L)+a_{p,1}(L)(\xi-1)+a_{p,2}(L)(\xi-1)^2+\cdots a_{p,p-1}(L)(\xi-1)^{p-1},
$$
for some integers $a_{p,n}(L)$'s. Then $(a_{p,n}(L) \bmod{p})\in \Z/p \Z \ (0\leq n \leq p-2)$ is uniquely 
determined by $\Phi_{\phi}(L)|_{t=\xi}$.
\end{lemma}
\begin{proof}
The ring $\Z[\xi]$ is isomorphic to $\Z[t]/(T(t))$ where
\begin{align*}
T(t)&=1+t+t^2+\cdots +t^{p-1}\\
&=\binom{p}{1}+\binom{p}{2}(t-1)+\cdots +\binom{p}{p}(t-1)^{p-1}.
\end{align*}
Since $p$ is an odd prime, $\binom{p}{k}$ is divisible by $p$ for $1\leq k \leq p-1$.
Hence $a_{p,k}$ is uniquely determined modulo $p$ by $\Phi_{\phi}(L)|_{t=\xi} \in \Z[\xi]$.
\end{proof}

Let $n$ be an integer, $Y$ be an indeterminate and
$f(Y)\in \Z/p\Z[Y,Y^{-1}]$ be a Laurent polynomial with $\Z/p\Z$ coefficients.
We see that
$$
\xi^{Y^n}=\exp (Y^n\log \xi)=\sum_{i=0}^{\infty}\frac{1}{i!}(Y^n \log \xi)^i\in \Q[\log \xi][[Y^n]].
$$
Hence,
$$
\xi^{f(Y)}\in \Q[\log \xi][\ldots, \xi^{Y^{-2}}, \xi^{Y^{-1}}, \xi^{Y^0}, \xi^{Y^1}, \xi^{Y^2}, \ldots].
$$
We prove that $\xi^{f(Y)}$ can be expressed as a sum of powers of $(\xi-1)$.

\begin{lemma}\label{lem4.2.}
Let $h_n$ be
$$
\sum_{i=1}^{n}(-1)^{i+1}\frac{(\xi-1)^i}{i}.
$$
We have that
$$
\xi^{f(Y)}\underset{(\xi-1)^{p-1}}\equiv
\sum_{i=0}^{p-2}\frac{h_{p-2}^i}{i!}{f(Y)}^{i}
\in \Z[\xi,\frac{1}{2},\frac{1}{3},\cdots ,\frac{1}{p-2},f(Y)].
$$
\end{lemma}
\begin{proof}
Since $h_{p-2}$ is congruent to $\log \xi$ modulo $(\xi-1)^{p-1}$,
we obtain the following formula: 
$$
\exp (h_{p-2} f(Y)) \equiv \xi^{f(Y)} \pmod{(\xi -1)^{p-1}}.
$$
Hence,
\begin{align*}
\xi^{f(Y)}&\equiv 1+h_{p-2} {f(Y)}+\frac{h_{p-2}^2}{2!}{f(Y)}^{2}+\cdots +\frac{h_{p-2}^{p-2}}{(p-2)!}{f(Y)}^{p-2}
\pmod{(\xi-1)^{p-1}}\\
&=\sum_{i=0}^{p-2}\frac{h_{p-2}^i}{i!}{f(Y)}^{i}
\in \Z[\xi,\frac{1}{2},\frac{1}{3},\cdots ,\frac{1}{p-2},f(Y)].
\end{align*}
\end{proof}

\section{Main Theorems}\label{203}

We prove the main Theorems. for any braids,
Vassiliev invariants derive from quandle shadow cosycle invariants 
by using $(\Z /p\Z)$-laurent polynomial quandles
and by the method similar to Theorem \ref{15}.

An {\it Alexander quandle} is defined to be $\Z / p \Z$ with a binary
operation given by 
$x \ast y=\om x+(1-\om)y$, 
where $\Z / p \Z$ means a cyclic group of order $p$ and 
$\om \in \Z / p \Z$ is neither $0$.
When $\om =1$, that is, $x\ast y=x$, we refer to the quandle as a {\it trivial quandle}.
For any $f(Y), \, g(Y) \in \Z/p\Z[Y,Y^{-1}]$, we define the following a binary operator:
$$
f(Y)\ast g(Y)=Y f(Y)+(1-Y) g(Y),
$$
we refer to the quandle as a {\it $(\Z /p\Z)$-laurent polynomial quandle}.

Moreover, T. Mochizuki obtained the following equation \cite{12:2003, 13:2005}:
$$
\phi(x,y,z)=(x-y)\cdot \Bigl(\frac{y^p-z^p-(y-z+\om^{-1}z)^p+(\om^{-1}z)^p}{p}\Bigr)
\in H^3_Q((\Z/p \Z,\om);\Z /p\Z).
$$
We call this cocycle the {\it Mochizuki $3$-cocycle}.
\begin{lemma}\label{201}
The following equation is a $3$-cocycle of a $(\Z /p\Z)$-laurent polynomial quandle:
$$
\phi'(f(Y),g(Y),h(Y))=(f(Y)-g(Y))(g(Y)-h(Y))(Y-1)\in \Z /p \Z[Y,Y^{-1}].
$$
\end{lemma}
\begin{proof}
We have that
\begin{align*}
\delta_3&(\phi')(f(Y),g(Y),h(Y),i(Y))=\Bigl((f(Y)-h(Y))(h(Y)-i(Y))\\
&-(Y f(Y)+(1-Y)g(Y)-h(Y))(h(Y)-i(Y))-(f(Y)-g(Y))(g(Y)-i(Y))\\
&+(Y f(Y)+(1-Y)h(Y)-(Yg(Y)+(1-Y) h(Y)))(Y g(Y)+(1-Y)h(Y)-i(Y))\\
&+(f(Y)-g(Y))(g(Y)-h(Y))-(Y f(Y)+(1-Y)i(Y)-(Yg(Y)+(1-Y) i(Y)))\\
& \ \ \times (Y g(Y)+(1-Y)i(Y)-(Yh(Y)+(1-Y) i(Y)))\Bigr)(Y-1) =0.
\end{align*}
Hence, $\phi'(f(Y),g(Y),h(Y)) \in \Z/p \Z[Y,Y^{-1}]$ is a $3$-cocycle.
\end{proof}

Let $b\in B_n$ be a braid,
let $\hat{D}$ be a braid diagram of the $n$-tops of the strands
and let $C'$ be a $\Z/p\Z$-coloring of $\mathcal{A}(\hat{D})\to \Z/p\Z$.
The closure of $b$ is a link $L$.
By definition (ii) of a quandle, there exists unique $z'\in \Z/p\Z$ in the unbounded region of 
a braid diagram of $L$.
As we define a quandle shadow cocycle invariant $\Phi_{\phi}(L)$,
we can define a quandle shadow cocycle invariant:
$$
\Phi_{\phi'}(b,Y)=\sum_{C'}\prod_{\x}W_{\phi '}(\x ;C')\in \Z[\Z/p \Z[Y,Y^{-1}])].
$$
However, in case of $\Phi_{\phi '}(b,Y)$,
the weighs of crossings (Figure \ref{10007}) must be calculated from the top to the bottom.
When we make closure, we regard the unbounded region of $L$ as $z'\in \Z/p\Z$.
Hence, $\Phi_{\phi'}(b, Y)$ is uniquely determined by $b$.
\begin{figure}[H]
 \centering
 \includegraphics[width=9cm]{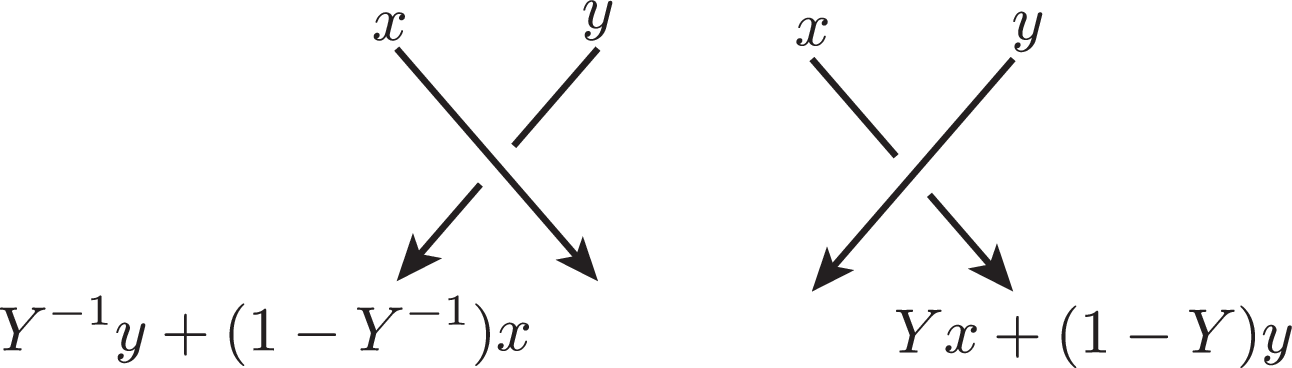}
 \caption{\rm The coloring at each crossing of braid $b$.} \label{10007}
\end{figure}


\begin{theorem}
Theorem \ref{100} is true:
let $p$ be an odd prime,
let $\xi={\rm exp}(2 \pi \sqrt{-1}/p)$,
let $\Z/p\Z \cong \langle t \mid t^p=1 \rangle$,
let $\phi'$ be a $3$-cocycle of a $(\Z /p\Z)$-laurent polynomial quandle $\Z/p \Z[Y,Y^{-1}]$
of Lemma \ref{201}
and let $\Phi_{\phi'}(b,Y)\in \Z[\Z/p\Z[Y,Y^{-1}]]$ be a quandle shadow cocycle invariant of a braid $b$
by using a coloring $C'$.
A quandle shadow cocycle invariant $\Phi_{\phi'}(b,Y)$ has the exponent a lauren polynomial
$\Z/p \Z[Y,Y^{-1}]$ of $t$.
Because of Lemma \ref{lem4.2.},
we have that
\begin{align*}
&\Phi_{\phi'}(b,Y)|_{t=\xi}\\
&\equiv \sum_{d=0}^{p-2}a_{p,d}(b,Y)(\xi-1)^d \pmod{(\xi-1)^{p-1}}\\
&\in \Z[\xi,\frac{1}{2},\frac{1}{3},\ldots ,\frac{1}{p-2},\phi']\\
&\cong \Z[t,\frac{1}{2},\frac{1}{3},\ldots ,\frac{1}{p-2}, \phi']/(1+t+t^2+\cdots +t^{p-1}).
\end{align*}
When we substitute $1/(1-\hbar)$ for $Y$ of $a_{p,d}(b,Y)\in \Z/p \Z[Y,Y^{-1}] \ (0\leq d \leq p-2)$,
the coefficient of $\hbar^k$ in
$
\displaystyle a_{p,d}(b,Y)|_{Y=1/(1-\hbar)}\in \Z/ p \Z[[\hbar ]]
$
is a $(\Z/p\Z)$-valued Vassiliev invariant of degree $k$ for any $0\leq d\leq p-2$.
\end{theorem}
\begin{proof}
Firstly, we show that 
$a_{p,d}(b,Y)|_{Y=1/(1-\hbar)} \in \Z/p \Z[[ \hbar ]]$
is a $(\Z/p \Z)$-valued Vassiliev invariants of degree $k$ of $b$.
Regard braids obtained by for the diagram $b_{p_1}$ with a double point $p_1$, replacing the double point with
a positive crossing and a negative crossing as $b_+$ and $b_-$ respectively.
When $\hbar=0$, we obtain the following equation:
$$
\Phi_{\phi'}(b_{p_1},1)=\Phi_{\phi'}(b_+,1)-\Phi_{\phi'}(b_-,1)= 0.
$$
Hence, $\displaystyle \Phi_{\phi'}(b_{p_1},Y)|_{Y=1/(1-\hbar)}$ is divisible by $\hbar$.
By using it repeatedly, if $b_{p_1p_2\cdots p_{k+1}}$ is a singular braid with exact $k+1$ singular points, then 
$\displaystyle \Phi_{\phi'}(b_{p_1p_2\cdots p_{k+1}},Y)|_{Y=1/(1-\hbar)}$
is divisible by $\hbar^{k+1}$.
Therefore, $a_{p,d}(b_{p_1p_2\cdots p_{k+1}},Y)|_{Y=1/(1-\hbar)}$ is divisible by $\hbar^{k+1}$ in
$\Z/p \Z[[\hbar]]$.
In particularly, the coefficient of $\hbar^k$ is equal to $0$ for such singular braids.

Secondly, because of Lemma \ref{13},
the value of $a_{p,d}(b,Y)\in \Z/p \Z$ is unique. 
\end{proof}

Let $L$ be a link.
The closure of $b=\sigma_{|x_1|}^{{\rm sign}(x_1)}\cdots \sigma_{|x_n|}^{{\rm sign}(x_n)}$ is a link the $L$.
Here, $(x_1,\ldots x_n)\in \Z^n$ is a minimal lattice presentation of $L$.
\begin{theorem}[{\cite{28:2010}}]
Let $\mathbb{L}$ be the set of all links modulo equivalence and
let $\mathbb{B}$ be the set of all braids modulo equivalence. 
There exists the following injective map $\sigma$ by using minimal lattice presentation of the $L$:
$$
\sigma:\mathbb{L}\hookrightarrow\mathbb{B}.
$$
\end{theorem}
We can calculate a quandle shadow cocycle invariant
$\Phi_{\phi'}(b,-1)=\Phi_{\phi'}(\sigma (L),-1)$ by using a $(\Z/p \Z)$-laurent polynomial quandle.
Moreover, we obtain Vassiliev invariants for a braid $\sigma (L)$, 
but this is not the Vassiliev invariants of the link $L$.

\begin{example}\label{110}
Let  $K$ be $(2,3)$-torus knot and $p=3$.
In addition, we consider the following quandle shadow cocycle invariant:
\begin{align*}
\Phi_{\phi''}(\sigma (K),Y)=&9+3t^{\frac{1}{3}(-Y+Y^2-Y^3)(Y-1)}
+2t^{\frac{1}{3}(-4Y+4Y^2-4Y^3)(Y-1)}+4t^{\frac{1}{3}(-1-Y^3)(Y-1)}\\
&+2t^{\frac{1}{3}(1-2Y+2Y^2-Y^3)(Y-1)}
+2t^{\frac{1}{3}(-4-4Y^3)(Y-1)}+2t^{\frac{1}{3}(-2+Y-Y^2-Y^3)(Y-1)}\\
&+t^{\frac{1}{3}(-2-2Y+2Y^2-4Y^3)(Y-1)}+t^{\frac{1}{3}(-Y+Y^2-Y^3)(Y-1)}\\
&+t^{\frac{1}{3}(-2-2Y+2Y^2-4Y^3)(Y-1)}\in \Z[t]/(t^3-1).
\end{align*}
Here, we use the following $3$-cocycle:
$$
\phi''(f(Y),g(Y),h(Y)):=\frac{1}{3}(f(Y)-g(Y))(g(Y)-h(Y))(Y-1).
$$
We have that
$$
\frac{\Phi_{\phi''}(K,-1)}{3}=5+2t+2t^2,
$$
and
\begin{align*}
\frac{a_{3,1}(\sigma (K),Y)}{3}&=(-2-2Y+2Y^2-4Y^3)(Y-1)\\
&\mapsto -6\hbar-16\hbar^2-36\hbar^3-70\hbar^4-122\hbar^5-196\hbar^6
+O(\hbar^7)\in \Z/3 \Z[[\hbar]].
\end{align*}
Because of Theorem \ref{100},
$-16\equiv 2 \in \Z/3 \Z$ is one of the $(\Z/3 \Z)$-valued Vassiliev invariants of degree $2$,
$-36 \equiv 0 \in \Z/3 \Z$ is one of the $(\Z/3 \Z)$-valued Vassiliev invariants of degree $3$,
$-70 \equiv 2 \in \Z/3 \Z$ is one of the $(\Z/3 \Z)$-valued Vassiliev invariants of degree $4$,
$-122 \equiv 1 \in \Z/3 \Z$ is one of the $(\Z/3 \Z)$-valued Vassiliev invariants of degree $5$ and
$-196 \equiv 2 \in \Z/3 \Z$ is one of the $(\Z/3 \Z)$-valued Vassiliev invariants of degree $6$.
\end{example}

We consider the following problem of surface $2$-knots:
\begin{problem}
$(1)$ Do singular curves corresponding to singular points exist ?\\
$(2)$ Can one define quantum invariants for 2-knots ?
\end{problem}


\section*{Acknowledgements}
The author expresses his sincere thanks to Special Appointment Professor 
T. Kanenobu (at Osaka Central Advanced Mathematical Institute) 
and Professor K. Tanaka (at Tokyo Gakugei University) for giving valuable advice.



\bigskip
\address{
SHINSHU UNIVERSITY \\
Faculty of Science \\
3-1-1 Asahi, Matsumoto, \\ 
Nagano 390-8621 \\
Japan
}
{k9250979@kadai.jp}


\end{document}